\documentclass[a4paper, oneside]{amsart}
\usepackage[utf8]{inputenc}

\usepackage{amsfonts}
\usepackage{amsmath}
\usepackage{amssymb}
\usepackage{amsthm}
\usepackage[a4paper]{geometry}
\usepackage{mathrsfs} 
\usepackage[dvipsnames]{xcolor}

\usepackage[pagebackref=true,colorlinks=true, linkcolor=MidnightBlue, citecolor=MidnightBlue]{hyperref}
\renewcommand*\backref[1]{\ifx#1\relax \else (Cited on #1) \fi}

\usepackage{appendix}
\usepackage{enumerate}
\usepackage{geometry}
\usepackage{tikz-cd} 

\theoremstyle{plain}
\newtheorem{theorem}{Theorem}
\newtheorem{proposition}[theorem]{Proposition}
\newtheorem{lemma}[theorem]{Lemma}
\newtheorem{corollary}[theorem]{Corollary}
\newtheorem{definition}[theorem]{Definition}
\newtheorem{remark}[theorem]{Remark}

\theoremstyle{definition}
\newtheorem{example}[theorem]{Example}

\numberwithin{theorem}{section}
\numberwithin{equation}{section} 

\newcommand{\set}[1]{\left\{ #1 \right\}}
\newcommand{\norm}[1]{\left \lVert  #1 \right \rVert}
\newcommand{\abs}[1]{\left\lvert #1 \right\rvert}

\newcommand{\domL}{\text{\normalfont dom}(\mathscr{L})}

\newcommand{\vertiii}[1]{{\left\vert\kern-0.25ex\left\vert\kern-0.25ex\left\vert #1 
    \right\vert\kern-0.25ex\right\vert\kern-0.25ex\right\vert}}

\newcommand{\Z}{\mathbb{Z}}

\newcommand{\R}{\mathbb{R}}
\newcommand{\N}{\mathbb{N}}

\newcommand{\calA}{\mathcal{A}}
\newcommand{\calB}{\mathcal{B}}

\newcommand{\calF}{\mathcal{F}}
\newcommand{\calG}{\mathcal{G}}

\newcommand{\calM}{\mathcal{M}}

\newcommand{\bbE}{\mathbb{E}}
\newcommand{\bbP}{\mathbb{P}}
\newcommand{\condexp}[2]{\bbE \left[ #1 \lvert #2 \right]}

\newcommand{\dom}[1]{\text{\upshape dom}(#1)}

\newcommand{\phientropy}{\mathbf{Ent}_\mu^\Phi}
\newcommand{\phidiv}{\text{\upshape div}^{\Phi}}

\title{Trajectorial dissipation of $\Phi$-entropies for interacting particle systems}
\author{Benedikt Jahnel}
\address{Institut f\"ur Mathematische Stochastik, Technische Universit\"at Braunschweig, Universit\"atsplatz 2,
38106 Braunschweig, Germany \& Weierstrass Institute for Applied Analysis and Stochastics\\
Mohrenstraße 39\\
10117 Berlin\\
Germany}
\email{benedikt.jahnel@tu-braunschweig.de}

\author{Jonas Köppl}
\address{Weierstrass Institute for Applied Analysis and Stochastics\\
Mohrenstraße 39\\
10117 Berlin\\
Germany}
\email{jonas.koeppl@wias-berlin.de}

\date{\today}

\begin{document}

\maketitle

\begin{abstract}
     A classical approach for the analysis of the longtime behavior of Markov processes is to consider suitable Lyapunov functionals like the variance or more generally $\Phi$-entropies. Via purely analytic arguments it can be shown that these functionals are indeed non-increasing in time under quite general assumptions on the process. 
     In this paper, we complement these classical results via a more probabilistic approach and show that dissipation is already present on the level of individual trajectories for spatially-extended systems of infinitely many interacting particles with arbitrary underlying geometry and compact local spin spaces.
    This extends previous results from the setting of finite-state Markov chains or diffusions in $\R^n$ to an infinite-dimensional setting with weak assumptions on the dynamics. 
\end{abstract}


\section{Introduction}\label{section:introduction}

There are many different techniques to study the long-time behavior of Markov processes that excel in different situations. One very common and powerful technique is the use of Lyapunov functionals, i.e., functionals that are monotone in time. An example of such a functional is the variance
\begin{align*}
    \mathbf{Var}_\mu(f) := \mathbb{E}_\mu[f^2]-\mathbb{E}_\mu[f]^2, \quad f \in L^2(\mu),
\end{align*}
where $\mu$ is an invariant measure for some Markov process with semigroup $(P_t)_{t \geq 0}$. If we now fix an observable $f$ and consider the function 
\begin{align*}
    [0,\infty) \ni t  \mapsto \mathbf{Var}_\mu(P_tf) \in [0,\infty),
\end{align*}
then it is easy to see that this is non-increasing and under some further assumptions one can even show that it is strictly decreasing for all non-constant observables $f$. This whole viewpoint, however, is purely based on functional analytic arguments and one does not even need to speak about the underlying Markov process itself to carry out the corresponding calculations. From a probabilistic point of view, this is somewhat unsatisfactory and we therefore want to specify this coarse and non-probabilistic approach with a finer, more probabilistic technique that allows us to obtain trajectorial versions of these results. Here, by trajectorial we mean results on the level of \textit{single realizations} of a stochastic process, as opposed to averaged quantities. Thereby our goal is to uncover more of the underlying probabilistic mechanisms behind the decay of variance, or more generally the decay of $\Phi$-entropies. 
For this, we will first briefly recall the notion of $\Phi$-entropies and then explain our main results and ideas with the help of the simple example of a continuous-time Markov chain on a finite state space. The rest of the article is then devoted to extending these ideas to the setting of spatially extended systems of infinitely many interacting particles as e.g.~considered in~\cite{liggett_interacting_2005}. 

\subsection{$\Phi$-entropies and their decay under Markovian dynamics}

Let $\Phi: I \to \R$ be a smooth and \textit{convex} function defined on a not necessarily bounded interval $I \subset \R$. Let $(E, \calB(E))$ be a Polish space equipped with its Borel $\sigma$-algebra and assume that $\mu$ is a probability measure on $(E, \calB(E))$. 
The $\Phi$-entropy functional is then defined on the set of $\mu$-integrable functions $f:E \to I$ by 
\begin{align*}
    \phientropy(f) 
    := 
    \int_E \Phi(f) d\mu - \Phi\left( \int_E f d \mu\right) 
    = 
    \mathbb{E}_\mu\left[\Phi(f)\right] - \Phi\left(\mathbb{E}_\mu\left[f\right]\right). 
\end{align*}
By Jensen's inequality one can immediately deduce that the $\Phi$-entropy functional takes its values in $\R_+ \cup \set{+\infty}$. 
Moreover, $\phientropy(f)$ vanishes if its argument is constant and if $\Phi$ is \textit{strictly} convex, then the converse is also true. 
For special choices of $\Phi$ one can recover the classical variance and relative entropy functionals since we have 
\begin{align*}
    \mathbf{Ent}_\mu^{u \mapsto u^2} = \mathbf{Var}_\mu, \quad \mathbf{Ent}_\mu^{u \mapsto u \log u} = h(\cdot \lvert \mu). 
\end{align*}

Now let $(X(t))_{t \geq 0}$ be a Markov process on our Polish space $E$ with associated semigroup $(P_t)_{t \geq 0}$ acting on $C_b(E;\R)$, the space of continuous and bounded real-valued functions on $E$. Let us assume that there exists an invariant probability measure $\mu$ and denote by $\mathscr{L}$ the generator of the semigroup $(P_t)_{t \geq 0}$ with domain $\dom{\mathscr{L}} \subset C_b(E;\R)$. 

By invariance of $\mu$ and Jensen's inequality one can now deduce that for all $f \in C_b(E;\R)$
\begin{align*}
    \phientropy(P_tf) 
    &= 
    \mathbb{E}_\mu\left[\Phi(P_tf)\right] - \Phi\left(\mathbb{E}_\mu\left[P_tf\right]\right)\leq 
    \mathbb{E}_\mu\left[P_t\Phi(f)\right] - \Phi\left(\mathbb{E}_\mu\left[f\right]\right)=
    \phientropy(f). 
\end{align*}
This tells us that the $\Phi$-entropy is non-increasing as a function of $t$ and can be used as a Lyapunov function. More precisely, with purely analytic arguments, one can even deduce the following general result about the decay of $\Phi$-entropies. 

\begin{proposition}[DeBruijn like property for Markov semigroups, \cite{chafai_entropies_2004}]\label{prop:debruijn-property}
Let $(X(t))_{t \geq 0}$ be a Markov process on a Polish space $E$ equipped with its Borel $\sigma$-algebra $\calB(E)$ and let $(P_t)_{t \geq 0}$ be the associated Markov semigroup with generator $\mathscr{L}$. 
Assume that $\mu$ is an invariant probability measure. Then, for any continuous and bounded function $f:E \to I$ and any $t>0$, it holds that \begin{align*}
    \partial_t \ \phientropy(P_tf) = \mathbb{E}_\mu \left[\Phi'(P_t f)\mathscr{L}(P_t f) \right] \leq 0. 
\end{align*}
\end{proposition}

This result is classical, but we nevertheless recall its short analytic proof. We will later also provide a more probabilistic argument to obtain the same result in the context of interacting particle systems. 

\begin{proof}
The chain rule and the definition of the generator $\mathscr{L}$ directly imply that 
\begin{align*}
    \partial_t \ \phientropy(P_tf) 
    = 
    \mathbb{E}_\mu \left[\Phi'(P_t f)\frac{d}{dt}(P_t f) \right] 
    = 
    \mathbb{E}_\mu \left[\Phi'(P_t f)\mathscr{L}(P_t f) \right].
\end{align*}
To see that the left-hand side is actually non-positive, it suffices to observe that the convexity of $\Phi$ implies via Jensen's inequality for conditional expectations
\begin{align*}
    \Phi(P_{t+s}g) \leq P_t(\Phi(P_sg))
\end{align*}
for any $s,t \geq 0$ and hence for all $f$ we have
\begin{align*}
    \phientropy(P_{t+s}f) \leq \phientropy(P_sf),
\end{align*}
by invariance of $\mu$. 
\end{proof}

By integrating one obtains the following classical corollary, which links exponential decay of $\Phi$-entropies and functional inequalities involving $\Phi$-entropies. 
\begin{corollary} 
In the setting of Proposition \ref{prop:debruijn-property}, the following two statements are equivalent. 
\begin{enumerate}[i.]
    \item There exists a constant $c>0$ such that for all $f\in\dom{\mathscr{L}}$
    \begin{align*}
        \phientropy(f) \leq -c \mathbb{E}_\mu \left[\Phi'(f)\mathscr{L}f \right].
    \end{align*}
    \item There exists a constant $c>0$ such that for all continuous and bounded $f:E \to I$
    \begin{align*}
        \phientropy(P_tf) \leq e^{-\frac{t}{c}}\phientropy(f).
    \end{align*}
\end{enumerate}
\end{corollary}
Note that, in the special case $\Phi: u \mapsto u^2$, one recovers the Poincaré inequality
\begin{align*}
    \mathbf{Var}_\mu(f) \leq -\frac{c}{2} \ \langle f, \mathscr{L}f \rangle_{L^2(\mu)},
\end{align*}
which is well-known to be equivalent to exponential $L^2$ ergodicity, see e.g.~\cite{holley_l2_1976}. For a more detailed review of $\Phi$-entropies and further results we refer the interested reader to \cite{chafai_entropies_2004}. 

\subsection{A finite state-space example for the trajectorial approach}\label{section:intro-trajectorial-idea}
As one can see, the results above can be obtained without even mentioning the underlying stochastic process and just dealing with the semigroup and its generator. We want to supplement this viewpoint with a more probabilistic approach that allows us to obtain a somewhat finer result on a trajectorial level, which also implies the classical results by taking expectations.

For simplicity, we will discuss the main ideas for the example of a continuous-time Markov chain on a finite state space. 
More precisely, let $(X_t)_{t\geq 0}$ be a Markov chain on a finite set $E$ with irreducible generator $\mathscr{L}$ and strictly positive invariant measure $\mu$. Hence, the corresponding Markov semigroup is given by the matrix exponential $(e^{t\mathscr{L}})_{t \geq 0}$.  Denote the underlying probability space by $(\Omega, \calA, \mathbb{P})$ and assume that $X_0 \sim \mu$ under $\mathbb{P}$. 



It is easy to check that for all bounded $f:[0,\infty) \times E \to \R$ such that for all $x \in E$ the partial derivatives $\partial_t f(\cdot, x)$ are continuous and bounded,  the process defined by 
\begin{align}\label{eqn:discrete-time-martingale-property}
    f(t, X_t) - \int_0^t (\partial_s + \mathscr{L})f(s,X_s)ds, \quad t \geq 0, 
\end{align}
is a martingale w.r.t.~the canonical filtration generated by $(X_t)_{t \geq 0}$, see e.g.~\cite[Lemma~IV.4.20]{rogers_diffusions_2000}. 

If we now fix a finite time horizon $T>0$ and consider the time-reversal $(\hat{X}_t)_{0 \leq t \leq T}$ of $(X_t)_{t\ge0}$, where $\hat{X}_t = X_{T-t}$, then under $\mathbb{P}$ the time-reversed process is again a time-homogeneous Markov process with generator $\hat{\mathscr{L}}$, where
\begin{align*}
    \hat{\mathscr{L}}(x,y) = \frac{\mu(y)}{\mu(x)}\mathscr{L}(y,x).
\end{align*}
A short calculation now shows that, for each bounded $g : E \to \R$ and $T>0$, the process $(P_{T-s}g(\hat{X}_s))_{0 \leq s\leq T}$ is a $((\hat{\calF}_t)_{0 \leq t \leq T}, \mathbb{P})$-martingale, where $\hat{\calF}_t = \sigma(X_{T-s}: \ 0 \leq s \leq t)$. Indeed, we can use the chain rule to calculate 
\begin{align*}
   \partial_t P_{T-t}g(x) = - \hat {\mathscr{L}} P_{T-t}g(x),
\end{align*}
so the correction term in (\ref{eqn:discrete-time-martingale-property}) vanishes.
Note that it is crucial to use the time-reversed process here, since the correction term does not cancel out if one uses the forward process. 

By convexity this directly implies that the time-reversed trajectorial $\Phi$-entropy, i.e., the process defined by 
\begin{align*}
    \Phi(P_{T-s}g(X(T-s))), \quad 0 \leq t \leq T,
\end{align*}
is a submartingale. This stochastic monotonicity  can be seen as a trajectorial version of the classical $\Phi$-entropy decay and after taking expectations with respect to $\mathbb{P}$, one obtains the classical results as in Proposition \ref{prop:debruijn-property}. Therefore, one can interpret these trajectorial results as the probabilistic version of the decay of $\Phi$-entropies. 

\medskip
The main work is now to establish that a similar argument can also be used to treat infinite-dimensional systems like the interacting particle systems we consider. 
To our best knowledge, the first results of this kind, in the context of diffusions in $\R^n$, have been achieved in \cite{fontbona_trajectorial_2016}. 
More recently, starting with \cite{karatzas_trajectorial_2020}, these results have been extended to more and more classes of Markov processes, including continuous time Markov chains on countable state spaces, see \cite{karatzas_trajectorial_2021}. The works \cite{kim_trajectorial_2022} and \cite{tschiderer_trajectorial_2023} are also in a similar spirit. 

\medskip
The setting will be made precise in Section \ref{section:setting-results}, but roughly speaking, we consider continuous-time Markov jump processes on general configuration spaces $\Omega = \Omega_0^S$, where $S$ is an arbitrary countable set and $\Omega_0$ is a compact Polish space. We will refer to the elements of $S$ as \textit{sites} and call $\Omega_0$ the \textit{local state-space}. 
In most examples considered in the literature, $S$ is the vertex set of some graph like the $d$-dimensional hypercubic lattice $\Z^d$, a tree or the Cayley graph of a group. This underlying spatial geometry dictates which particles can interact with each other and we are therefore not in the setting of mean-field systems but in an infinite-dimensional setting. This of course brings with it its own set of technical difficulties which need to be dealt with for making the time-reversal arguments work. 

The main technical difficulties come from making sure that the time-reversal is again a well-defined interacting particle system and from obtaining a description of its generator. This is made possible by assuming some local regularity of the local conditional distributions of the time-stationary measure $\mu$. Namely, by the assumption that $\mu$ is actually a Gibbs measure with respect to a quasilocal specification that additionaly satisfies a certain smoothness condition. This condition is e.g.~satisfied if the specification is given in terms of a potential $\Phi = (\Phi_B)_{B \Subset S}$ such that 
\begin{align*}
    \sup_{x \in S}\sum_{B \Subset S: \ B \ni x}\abs{B}\norm{\Phi_B}_\infty < \infty. 
\end{align*}
Note that this condition is for example satisfied for any translation-invariant finite-range potential, so our theory applies to a fairly large class of models. 

\subsection{Organization of the manuscript}
The rest of this article is organized as follows. We will first collect the necessary notation and formulate our main results in Section \ref{section:setting-results}.
Then, as a first step, we investigate the time-reversal of interacting particle systems in equilibrium in Section \ref{section:time-reversal-ips} with the main goal of obtaining an explicit representation of the (formal) generator of the time-reversed dynamics. 
In Section \ref{section:pathwise-phi-entropy-decay}, we will then apply these results to establish pathwise properties of general $\Phi$-entropy functionals and derive the classical DeBrujin-like decay property as a corollary. 

\section{Setting and main results}\label{section:setting-results}

Let $(\Omega_0, \calB_0)$ be a compact Polish space equipped with its Borel $\sigma$-algebra and $\lambda_0$ a probability measure on $(\Omega_0, \calB_0)$, which will serve as our reference measure.
We will consider Markovian dynamics on the configuration space $\Omega = \Omega_0^{S}$, where $S$ is some countable set whose elements we will refer to as \textit{sites}. 
In most applications this will be the set of vertices of some graph, e.g.~$\Z^d$ or a tree. 
We equip $\Omega$ with the product topology and corresponding Borel $\sigma$-algebra $\calF$. Note that $\calF$ coincides with the product $\sigma$-algebra $\otimes_{x \in S}\calB_0$.
For $\Delta \subset S$ we will also write $\Omega_\Delta := \Omega_0^\Delta$ for the set of partial configurations. We will also equip $\Omega_\Delta$ with the product $\sigma$-algebra and the probability measure $\lambda_\Delta = \otimes_{x \in \Delta}\lambda_0$. 
For $\Lambda \subset S$, let $\calF_\Lambda$ be the sub-$\sigma$-algebra of $\calF$ that is generated by the projections $\omega \mapsto \omega_\Delta \in \Omega_\Delta$ for $\Delta \Subset \Lambda$, where we write $\Subset$ to signify that a set is a \textit{finite} subset of another set. 
For $\Delta \subset S$ and (partial) configurations $\eta_{\Delta^c} \in \Omega_{\Delta^c}$ and $\xi_\Delta \in \Omega_\Delta$, we will write $\xi_\Delta\eta_{\Delta^c}$ for the configuration that is defined on all of $S$ and agrees with $\eta_{\Delta^c}$ on $\Delta^c$ and with $\xi_\Delta$ on $\Delta$. 
For a topological space $E$, we will denote its Borel $\sigma$-algebra by $\calB(E)$ and the space of continuous real-valued functions on $E$ by $C(E)$. The space of non-negative measures on $E$, or more precisely on $\calB(E)$,  will be denoted by $\calM(E)$ and is equipped with the topology of weak-convergence. The total variation distance on $\calM(E)$ will be denoted by $\norm{\cdot}_{\text{TV}}$. 


\subsection{Interacting particle systems and Gibbs measures}
\subsubsection{Interacting particle systems}
We will consider time-continuous Markovian dynamics on $\Omega$, namely interacting particle systems characterized by time-homogeneous generators $\mathscr{L}$ with domain $\text{dom}(\mathscr{L}) \subset C(\Omega)$ and the associated Markovian semigroup $(P_t)_{t \geq 0}$ on $C(\Omega)$. 
For interacting particle systems we adopt the notation and exposition of the standard reference \cite[Chapter~I]{liggett_interacting_2005}. 

In our setting, the generator $\mathscr{L}$ is given by a collection of transition measures $(c_\Delta(\cdot, d\xi))_{\Delta \Subset S}$ in finite volumes $\Delta \Subset S$, i.e., mappings 
\begin{align*}
    \Omega \ni \eta \mapsto c_\Delta(\eta, d\xi_\Delta) \in \calM(\Omega_\Delta). 
\end{align*}
These transition measures can be interpreted as the infinitesimal rates at which the particles inside $\Delta$ switch from the configuration $\eta_\Delta$ to $\xi_\Delta$, given that the rest of the system is currently in state $\eta_{\Delta^c}$. The full dynamics of the interacting particle system is then given as the superposition of these local dynamics,
\begin{align}\label{eqn:formal-generator}
    \mathscr{L}f(\eta) = \sum_{\Delta \Subset S}\int_{\Omega_\Delta}\left[f(\xi_\Delta \eta_{\Delta^c})-f(\eta)\right]c_\Delta(\eta, d\xi_\Delta). 
\end{align}
In \cite[Chapter~I]{liggett_interacting_2005} it is shown that the following conditions are sufficient to guarantee well-definedness. 
\begin{enumerate}[\bfseries (L1)]
    \item For each $\Delta \Subset \Omega$ the mapping 
    \begin{align*}
        \Omega \ni \eta \mapsto c_\Delta(\eta, d\xi_\Delta) \in \calM(\Omega_\Delta)
    \end{align*}
    is continuous. 
    \item The total rate at which a single particle switches its state is uniformly bounded, i.e.,
    \begin{align*}
        \sup_{x \in S}\sum_{\Delta \ni x}\sup_{\eta \in \Omega}c_\Delta(\eta, \Omega_\Delta) < \infty. 
    \end{align*}
    \item The total influence of all other particles on the transition rates of a single particle is uniformly bounded, i.e.,
    \begin{align*}
        M := \sup_{x \in S}\sum_{\Delta \ni x}\sum_{y \neq x}\delta_y c_\Delta < \infty,
    \end{align*}
    where 
    \begin{align*}
       \delta_y c_\Delta := \sup\set{\norm{c_\Delta(\eta, \cdot) - c_\Delta(\xi, \cdot)}_{\text{TV}}: \ \eta_{y^c} = \xi_{y^c}}.
    \end{align*}
\end{enumerate}

Under these conditions, the core of the operator $\mathscr{L}$ is given by 
\begin{align*}
    D(\Omega) := \set{f \in C(\Omega): \ \vertiii{f} := \sum_{x \in S}\delta_x(f) < \infty},
\end{align*}
where for $x \in S$
\begin{align*}
    \delta_x(f) := \sup_{\eta, \xi: \ \eta_{x^c} = \xi_{x^c}}\abs{f(\eta) - f(\xi)}
\end{align*}
is the oscillation of a function $f:\Omega \to \R$ at the site $x$. In addition, one can show the following estimates for $\mathscr{L}$ and the action of the semigroup $(P_t)_{t \geq 0}$ generated by $\mathscr{L}$. We will need these later on. 

\begin{lemma}\label{lemma:growth-estimate-triple-norm}
Assume that the generator $\mathscr{L}$ satisfies $\mathbf{(L1)}-\mathbf{(L3)}$ and denote by $(P_t)_{t \geq 0}$ the semigroup generated by $\mathscr{L}$. 
\begin{enumerate}[i.]
    \item For $f\in D(\Omega)$ we have $P_tf \in D(\Omega)$ for all $t \geq 0$ and
    \begin{align*}
        \vertiii{P_tf} \leq \exp\left((M-\varepsilon)t\right)\vertiii{f}.
    \end{align*}
    \item For all $f \in D(\Omega)$ it holds that 
    \begin{align*}
        \norm{\mathscr{L}f}_\infty 
        \leq 
        \left(
        \sup_{x \in S}\sum_{\Delta \ni x}\sup_{\eta \in \Omega}c_\Delta(\eta, \Omega_\Delta) 
        \right) \vertiii{f}. 
    \end{align*}
\end{enumerate}
The constants are explicitly given by 
\begin{align*}
    M 
    &=
    \sup_{x \in S} \sum_{\Delta \ni x}\sum_{y \neq x}\delta_y c_{\Delta} < \infty,
    \\\
    \varepsilon 
    &= 
    \inf_{x \in S} \inf_{\eta,\zeta\colon \eta_{x^c} = \zeta_{x^c}, \eta_x \neq \zeta_x} 
        \sum_{\Delta \ni x} \left(
            \sum_{\xi_{\Delta}: \xi_x = \zeta_x} c_{\Delta}(\eta, \xi_{\Delta})
            +
            \sum_{\xi_{\Delta}: \xi_x = \eta_x} c_{\Delta}(\zeta, \xi_{\Delta})
        \right). 
\end{align*}
\end{lemma}

\begin{proof}
    Combine the results from Proposition 3.2(a) and Theorem 3.9.(d) in \cite[Chapter I]{liggett_interacting_2005}. 
\end{proof}

For our purposes, the mere well-definedness of an interacting particle system is not sufficient and we need to assume some more regularity. All the additional assumptions we put will be used to make the generator of the time-reversal well-defined. 
\begin{enumerate}[\bfseries (R1)]
    \item For each $\Delta$ and $\eta \in \Omega$ the measure $c_\Delta(\eta, d\xi_\Delta)$ is absolutely continuous w.r.t.~ the reference measure $\lambda_\Delta(d\xi_\Delta)$ with density $c_\Delta(\eta, \cdot)$.  
    \item For each $\Delta \in \Omega$ the map 
    \begin{align*}
        \Omega \times \Omega_\Delta \ni (\eta, \xi_\Delta) \mapsto c_\Delta(\eta, \xi_\Delta) \in \R
    \end{align*}
    is continuous w.r.t.~the product topology. 
    \item The total rate of transition for a single site is uniformly bounded from above 
    \begin{align*}
        \sup_{x \in S}\sum_{\Delta \ni x}\sup_{\eta \in \Omega}\norm{c_\Delta(\eta, \cdot)}_\infty < \infty. 
    \end{align*}
    \item 
    The condition $\mathbf{(L3)}$ is satisfied, i.e.,
    \begin{align*}
        \sup_{x \in S}\sum_{\Delta \ni x}\sum_{y \neq x} \delta_yc_\Delta < \infty.
    \end{align*}
    \item There exists an $R>0$ such that for all $\Delta \Subset \Z^d$ with $\abs{\Delta} > R$ we have 
    \begin{align*}
        \sup_{\eta \in \Omega,\xi_\Delta \in \Omega_\Delta}c_\Delta(\eta,\xi_\Delta) = 0. 
    \end{align*}
\end{enumerate}

We will comment on where and why we need these assumptions and their connection to the classical conditions $(\mathbf{L1})-(\mathbf{L3})$ at the end of Section \ref{subsection:def-specification}, after we have stated our assumptions on the local conditional distribution of the time-stationary measure $\mu$. 

\subsubsection{Gibbs measures and the DLR formalism}\label{subsection:def-specification}
We will mainly be interested in the situation where the process generated by $\mathscr{L}$ admits a time-stationary measure $\mu$ with a well-behaved local representation, namely that $\mu$ is a Gibbs measure w.r.t.~to a sufficiently nice specification $\gamma$. Let us  therefore first recall the general definition of a specification.  

\begin{definition}
A specification $\gamma = (\gamma_\Lambda)_{\Lambda \Subset S}$ is a family of probability kernels $\gamma_{\Lambda}$ from $\Omega _{\Lambda^c}$ to $\calM_1(\Omega)$ that additionally satisfies the following properties. 
\begin{enumerate}[i.]
    \item Each $\gamma_{\Lambda}$ is \emph{proper}, i.e., for all $B \in \calF_{\Lambda^c}$ it holds that 
    \begin{align*}
        \gamma_\Lambda(B | \cdot) = \mathbf{1}_B(\cdot). 
    \end{align*}
    \item The probability kernels are \emph{consistent} in the sense that if $\Delta \subset \Lambda \Subset S$, then for all $A \in \calF$
    \begin{align*}
        \gamma_{\Lambda}\gamma_{\Delta}(A|\cdot)
        = 
        \gamma_\Lambda(A|\cdot),
\end{align*}
where the concatenation of two probability kernels is defined as usual via 
\begin{align*}
    \gamma_\Lambda \gamma_\Delta(A \lvert \eta) = \int_\Omega \gamma_\Delta(A \lvert \omega) \gamma_\Lambda(d\omega \lvert \eta). 
\end{align*}
\end{enumerate}
\end{definition}

For the existence and further properties of Gibbs measures with specification $\gamma$ one needs to impose some conditions on the specification $\gamma$. One sufficient condition for the existence of a Gibbs measure for a specification $\gamma$ is \textit{quasilocality}, see e.g.~\cite{friedli_statistical_2017} or \cite{georgii_gibbs_2011}. 
For the following sections we will need to assume some more regularity for the specification $\gamma$. In particular, these assumptions will guarantee that $\gamma$ is quasilocal. 

\begin{enumerate}[\bfseries (S1)]
    \item For each $\Delta \Subset S$ and $\eta \in \Omega$, the probability measure $\gamma_\Delta(d\xi_\Delta \lvert \eta)$ is absolutely continuous w.r.t.~the reference measure $\lambda_\Delta(d\xi_\Delta)$ with density $\gamma_\Delta(\cdot \lvert \eta)$. 
    \item For all $\Delta \Subset S$, the map 
    \begin{align*}
        \Omega \times \Omega_\Delta \ni (\eta, \xi_\Delta) \mapsto \gamma_\Delta(\xi_\Delta \lvert \eta_{\Delta^c}) \in [0, \infty)     \end{align*}
    is continuous (w.r.t.~the product topology). 
    \item The conditional densities on the single spin spaces are uniformly bounded away from zero and infinity, i.e.,
    \begin{align*}
            0 < \delta \leq \inf_{x \in S}\inf_{\eta \in \Omega}\gamma_{x}(\eta_x|\eta_{x^c})  \leq \sup_{x \in S}\sup_{\eta \in \Omega}\gamma_x(\eta_x \lvert \eta_{x^c}) \leq \delta^{-1} <\infty. 
    \end{align*}
    \item We have 
    \begin{align*}
        \sup_{x \in S}\sum_{\Delta \ni x:\ c_\Delta >0} \sum_{y \neq x}\delta_y\gamma_\Delta < \infty,
    \end{align*}
    where 
    \begin{align*}
        \delta_y\gamma_\Delta = \sup\set{\norm{\gamma_\Delta(d\xi_\Delta \lvert \eta) - \gamma_\Delta(d\xi_\Delta \lvert \zeta)}_{\text{TV}}: \ \eta_{y^c}=\zeta_{y^c} }. 
    \end{align*}
\end{enumerate}

\begin{remark}
    Now that we have stated all of the conditions that we need, let us briefly comment on why and where we need them. 
    \begin{enumerate}[i.]
        \item Assumption $\mathbf{(R3)}$ clearly implies $\mathbf{(L2)}$ and together with $\mathbf{(R4)}$ ensures that the interacting particle system is well-defined. 
        \item Assumption $\mathbf{(R1)}$ and $\mathbf{(S1)}$ allow us to write down the local transition density of the time-reversal and $\mathbf{(S3)}$ makes sure that we are not performing a division by zero. 
        \item The further regularity assumptions $\mathbf{(R3)}$, $\mathbf{(R5)}$ $\mathbf{(S3)}$, $\mathbf{(S4)}$ and the continuity assumptions $\mathbf{(R2)}$ and $\mathbf{(S2)}$ make sure that the local transition densities of the time-reversal also give rise to a well-defined interacting particle system. 
        \item The quantity in $\mathbf{(S4)}$ is similar to the classical Dobrushin uniqueness condition, see \cite{georgii_gibbs_2011}. However, we only need it to be finite and not strictly smaller than one. 
    \end{enumerate}
\end{remark}

\begin{example}
    One particular class of models to which our theory can be applied to are spin systems for which the specification $\gamma$ is defined via a  potential $\Phi = (\Phi_B)_{B \Subset S}$ that satisfies
    \begin{align*}
        \sup_{x \in S}\sum_{B \ni x} \abs{B}\norm{\Phi_B}_\infty < \infty, 
    \end{align*}
    and where the rates are of the form 
    \begin{align*}
        c_\Delta(\eta, \xi_\Delta) 
        =
        \begin{cases}
            \exp\left(-\beta \sum_{B: B \cap \Delta \neq \emptyset}\Phi_B(\xi_\Delta \eta_{\Delta^c})\right), \quad &\text{if } \abs{\Delta} = 1,
            \\\
            0, &\text{otherwise}. 
        \end{cases}
    \end{align*}
    Instead of these single-site updates one could also consider updates in larger regions with a bounded diameter. 
    Then the rates satisfy $\mathbf{(R1)}-\mathbf{(R5)}$ and the specification satisfies $\mathbf{(S1)}-\mathbf{(S4)}$, as one can see by using similar arguments as in the proof of \cite[Lemma 6.28]{friedli_statistical_2017}. 
\end{example}

\subsection{The time-reversal of an interacting particle system}

In the notation of above, assume that $\mu \in \mathscr{G}(\gamma)$ is Gibbs measure for a quasilocal specification $\gamma$, i.e., assume that $\mu$ satisfies the DLR equations \begin{align*}
    \mu(f)=\mu(\gamma_{\Lambda}(f| \cdot))
\end{align*} 
for all $\Lambda\Subset S$ and bounded measurable functions $f$. 
Further assume that $\mu$ is time-stationary with respect to the Markovian dynamics with generator $\mathscr{L}$.  
Denote the semigroup generated by $\mathscr{L}$ by $(P_t)_{t \geq 0}$ and the corresponding process on $\Omega$ by $(\eta(t))_{t \geq 0}$. As discussed in Section \ref{section:intro-trajectorial-idea}, for each fixed $T>0$ the process $(\eta(T-t))_{0 \leq t \leq T}$ is again a time-homogeneous Markov process and under some suitable assumptions its associated semigroup has a generator $\hat{\mathscr{L}}$. But what does this generator look like? 
For general Markov processes it is not possible to give a closed form expression, but in our case we can use the special structure of  $\mathscr{L}$ as the superposition of local dynamics in finite volumes. In each of these finite volumes, it is clear how the time-reversal w.r.t.\ $\mu$ should look and we can hope that we can again write $\hat{\mathscr{L}}$ as the superposition of finite volume processes. 
With this Ansatz, the probabilistic intuition dictates the educated guess 
\begin{align}\label{eqn:time-reversal-densities}
    \hat{c}_\Delta(\eta, \xi_\Delta) = c_\Delta(\xi_\Delta \eta_{\Delta^c}, \eta_\Delta) \frac{\gamma_\Delta(\xi_\Delta  \lvert \eta_{\Delta^c})}{\gamma_\Delta(\eta_\Delta \lvert \eta_{\Delta^c})} 
\end{align}
for the transition densities appearing in the generator of the time-reversed interacting particle system.
However, at this stage, it is not obvious that the generator of the time-reversed system is again of the form (\ref{eqn:formal-generator}) and has precisely these rates. For Markov processes on finite state spaces this is an easy calculation but we have to put in some more work, which will be carried out in Section \ref{section:time-reversal-ips}. The main result there is the following. 

\begin{proposition}[Time-reversal generator]\label{proposition:time-reversal-ips}
   Let the rates of an interacting particle system with generator $\mathscr{L}$ satisfy $\mathbf{(R1)}-\mathbf{(R5)}$ and assume that $\mu$ is a time-stationary measure for the corresponding Markov semigroup $(P(t))_{t \geq 0}$ on $C(\Omega)$ that is generated by $\mathscr{L}$ such that $\mu$ is a Gibbs measure w.r.t.~a specification $\gamma$ that satisfies $\mathbf{(S1)}-\mathbf{(S4)}$.
   Then, the time-reversed process has a generator $\hat{\mathscr{L}}$ whose transition densities (w.r.t.~the reference measure $\lambda_0$) are given by 
\begin{align*}
    \hat{c}_\Delta(\eta, \xi_\Delta) = c_\Delta(\xi_\Delta \eta_{\Delta^c}, \eta_\Delta) \frac{\gamma_\Delta(\xi_\Delta  \lvert \eta_{\Delta^c})}{\gamma_\Delta(\eta_\Delta \lvert \eta_{\Delta^c})}. 
\end{align*}
\end{proposition}
The proof of this can be found at the end of Section \ref{section:time-reversal-ips}. 

\subsection{Trajectorial decay of $\Phi$-entropies}

With this auxiliary result at hand, we can then obtain the following result, which describes the dissipation of general $\Phi$-entropies on a trajectorial level. Before we state the theorem, let us  introduce some further notation to express the main equation in a cleaner way.
The Bregman $\Phi$-divergence associated with $\Phi : I \to \R$ is defined as 
\begin{align*}
    \phidiv (p \lvert q) := \Phi(p) - \Phi(q) - (p-q)\Phi'(q), \quad p,q \in I. 
\end{align*}
This is precisely the difference between the value of $\Phi$ at the point $p$ and the value of the first-order Taylor expansion of $\Phi$ around $q$, evaluated at $p$ and is non-negative, since we assumed that $\Phi$ is convex. Bregman divergences are sometimes also referred to as Bregman distances, despite not being a metric since they are in general not symmetric and do not satisfy the triangle inequality. They are however still useful for applying techniques from optimization theory in more general contexts, e.g.~in statistical learning theory \cite{banerjee_clustering_2005}. 

Note that we now have to be careful with the probability space and filtration we are working with, since we are talking about results on a trajectorial level. 
\begin{theorem}[Trajectorial decay of $\Phi$-entropies]\label{thm:trajectorial-decay-phi-entropy}
    Let $(\Omega, \calA, \mathbb{P})$ be a probability space on which the interacting particle system $(\eta(s))_{s \geq 0}$ is defined.
    Denote the generator of the interacting particle system by $\mathscr{L}$ and assume that its rates satisfy $\mathbf{(R1)}-\mathbf{(R5)}$ and assume that under $\mathbb{P}$ we have $\eta_0 \sim \mu$, where $\mu$ is a time-stationary measure for the corresponding Markov semigroup $(P_t)_{t \geq 0}$ on $C(\Omega)$ that is generated by $\mathscr{L}$ such that $\mu$ is a Gibbs measure w.r.t.~a specification $\gamma$ that satisfies $\mathbf{(S1)}-\mathbf{(S4)}$.
    Then, for any $f \in D(\Omega)$ and $T>0$, the process defined by 
    \begin{align}\label{thm:definition-phi-likelihood-process}
        L^{\Phi, f}(s) := \Phi(P_{T-s}f(\eta_{T-s})), \quad 0 \leq s \leq T,
    \end{align}
    is a $((\hat{\calG}_t)_{0 \leq t \leq T},\mathbb{P})$-submartingale, where $\hat{\calG}_t = \sigma(\eta(T-s): \ 0 \leq s \leq t)$. Its Doob--Meyer decomposition is given by 
    \begin{align}\label{thm:doob-meyer-decomposition}
       &L^{\Phi, f}(t) = M^{\Phi, f}(t)  
       \\\
       &\qquad+\int_0^t \sum_{\Delta \Subset S}\int_{\Omega_\Delta}\hat{c}_\Delta(\eta(T-s), \xi_\Delta)\phidiv ( P_{T-s}(f(\xi_\Delta \eta_{\Delta^c}(T-s))\lvert P_{T-s}f(\eta(T-s))) \lambda_\Delta(d\xi_\Delta)  ds. 
       \nonumber
    \end{align}
\end{theorem}
The proof of this theorem can be found at the end of Section \ref{section:pathwise-phi-entropy-decay}. As mentioned before, by taking expectations one recovers the classical DeBrujin-like decay of $\Phi$-entropies as stated in Proposition \ref{prop:debruijn-property}.

For the sake of concreteness, let us write out the result  explicitly for one of the simplest cases, namely the trajectorial decay of variance, corresponding to $\Phi: u \mapsto u^2$. 

\begin{corollary}[Trajectorial decay of variance]
In the setting of Theorem \ref{thm:trajectorial-decay-phi-entropy}, we have that, for any $f \in D(\Omega)$ and $T>0$, the process defined by $(P_{T-s}f(\eta_{T-s}))^2$, $0 \leq s \leq T$, is a $((\hat{\calG}_t)_{0 \leq t \leq T},\mathbb{P})$-submartingale, where $\hat{\calG}_t = \sigma(\eta(T-s): \ 0 \leq s \leq t)$. Its Doob--Meyer decomposition is given by 
    \begin{align*}
        (P_{T-s}f(\eta_{T-s}))^2
        =
        M^f(t) 
        + 
        \int_0^t \sum_{\Delta \Subset S} \int_{\Omega_\Delta}\hat{c}_\Delta(\eta(T-s), \xi_\Delta)\left[f(\xi_\Delta\eta_{\Delta^c}(T-s)) - f(\eta(T-s))\right]^2\lambda_\Delta(d\xi_\Delta) ds. 
    \end{align*}
\end{corollary}

\subsection{Outlook}
Even though we were able to show the trajectorial decay for the relative entropy under quite general assumptions on the dynamics, these results are not fully satisfactory in the context of statistical mechanics. The usually more interesting Lyapunov functional in this setting is the so-called \textit{relative entropy density}, as e.g.~ considered in \cite{friedli_statistical_2017}, which is not only defined for measures $\nu$ that are absolutely continuous w.r.t.~ $\mu$. 
Therefore, it would be much more natural to work with this functional $h(\cdot \lvert \mu):\calM_1^{inv}(\Omega) \to \R$ instead 
and one can show that is also a Lyapunov function for interacting particle systems under quite general assumptions, see \cite{jahnel_dynamical_2022}, 
but it is somewhat unclear how to even formulate conjectures about the trajectorial properties of this functional, since one cannot naively evaluate it pointwise. 

As we already saw in the case of a continuous time Markov chain on a finite state space, the main ingredient for this type of result is to obtain an explicit description of the generator of the time-reversed process. Another class of processes that could be of interest and is not covered by our results are systems which evolve continuously on their single spin spaces, as opposed to our pure-jump processes. The first example that comes to mind are of course systems of (infinitely-many) interacting diffusions, e.g.~ indexed by $\Z^d$. We expect that, if a given system of interacting diffusions admits a Gibbs measure as an invariant probability measure, then a combination of the techniques in \cite{karatzas_trajectorial_2020} and this article should yield analogous results – of course under some suitable regularity conditions on the coefficients.

\section{The time-reversed interacting particle system and its generator}\label{section:time-reversal-ips}

The main goal of this section is to prove Proposition \ref{proposition:time-reversal-ips}, thereby establishing that the generator of the time-reversal is indeed given by $\hat{\mathscr{L}}$. For this we will need to establish some regularity properties for the transition densities as defined in (\ref{eqn:time-reversal-densities}).

\subsection{Upper and lower bounds for the conditional densities}
Since we will need to deal with quotients involving the conditional densities $\gamma_\Delta$ on arbitrary finite subsets $\Delta \Subset S$, we will need to lift the upper and lower bounds from $\mathbf{(S3)}$ to this more general case. This is essentially the content of the following lemma. 

\begin{lemma}
    Let $\gamma$ be a specification that satisfies $\mathbf{(S1)}$ and $\mathbf{(S3)}$.
    Then, there exists a constant $C \in (0,\infty)$ such that for all $\Delta \Subset S$ we have the estimate
    \begin{align*}
        e^{-C\abs{\Delta}} 
        \leq  
        \inf_{ \eta \in \Omega}\gamma_\Delta(\eta_\Delta \lvert \eta_{\Delta^c}) 
        \leq 
        \sup_{\eta \in \Omega} \gamma_\Delta(\eta_\Delta \lvert \eta_{\Delta^c}) 
        \leq 
        e^{C\abs{\Delta}}. 
    \end{align*}
    This constant is precisely given by $C = \abs{\log \delta}$. 
\end{lemma}

\begin{proof}
    For this, fix an enumeration $i_1, \dots, i_k$ of the elements of $\Delta$ and introduce the notation 
    \begin{align*}
        [i_j, i_k] := \set{i_j, i_{j+1},\dots, i_k}, \quad 1 \leq j \leq k. 
    \end{align*}
    With this notation at hand, we can use the chain rule for conditional probability densities to write
    \begin{align}\label{proof:chain-rule}
        \gamma_\Delta(\eta_{\Delta}\lvert \eta_{\Delta^c})
        =
        \prod_{j=1}^k\gamma_{[i_1,i_j]}(\eta_{i_j}\lvert\eta_{[i_{j+1},i_k]}\eta_{\Delta^c}),
    \end{align}
    where $\gamma_{[i_1,i_j]}(\eta_{i_j}\lvert\eta_{[i_{j+1},i_k]}\eta_{\Delta^c})$ is the marginal conditional density of the measure $\gamma_{[i_1,i_j]}(d \eta_{[i_1,i_j]} \lvert\eta_{[i_{j+1},i_k]} \eta_{\Delta^c})$ w.r.t.~the site $i_j$. 
    But, using consistency of the specification $\gamma$, we have 
     \begin{align*}
        \gamma_{[i_1,i_j]}(\eta_{i_j}\lvert\eta_{[i_{j+1},i_k]}\eta_{\Delta^c})
        &=\int\gamma_{[i_1,i_j]}(d\xi_{[i_1,i_j]}\lvert\eta_{[i_{j+1},i_k]}\eta_{\Delta^c})\gamma_{i_j}(\eta_{i_j}\lvert\xi_{[i_1,i_{j-1}]}\eta_{[i_{j+1},i_k]}\eta_{\Delta^c}),
    \end{align*}
   which is, by assumption, upper bounded by $\delta^{-1}$ and lower bounded by $\delta$.
    In conjunction with the representation (\ref{proof:chain-rule}) this implies the desired upper and lower bound where the constant $C$ is explicitly given by $C= \abs{\log(\delta)}$. 
\end{proof}

As a corollary we now get the following estimate for the quotients that appear in the definition of transition density of the time-reversal (\ref{eqn:time-reversal-densities}).

\begin{lemma}\label{lemma:non-nullness-inequality}
    Let $\Delta \Subset S$ and $\gamma$ be a specification that satisfies $\mathbf{(S1)}$ and $\mathbf{(S3)}$. Then, for all $\Delta \Subset S$, $\eta \in \Omega$ and $\xi_\Delta \in \Omega_\Delta$, we have 
    \begin{align*}
        0
        <
        e^{-2C \abs{\Delta}} 
        \leq 
        \frac{\gamma_\Delta(\xi_\Delta \lvert \eta_{\Delta^c})}{\gamma_\Delta(\eta_\Delta \lvert \eta_{\Delta^c})} 
        \leq 
        e^{2C \abs{\Delta}}.
    \end{align*}
\end{lemma}

\subsection{The switching lemma}
Now that we can be sure that the densities as in (\ref{eqn:time-reversal-densities}) are actually well-defined and we are not performing a division by zero, we can start showing that $\hat{\mathscr{L}}$ is indeed the generator of the time-reversed process. The main technical tool will be the following lemma. 

\begin{lemma}\label{generalized-switching-lemma}
Let the rates of a well-defined interacting particle system with generator $\mathscr{L}$ satisfy $\mathbf{(R1)}$ and assume that $\mu$ is a time-stationary measure for $\mathscr{L}$ and $\mu$ is a Gibbs measure w.r.t.~a specification $\gamma$ that satisfies $\mathbf{(S1)}$ and $\mathbf{(S3)}$.
Then, we have for all bounded and measurable $f,g:\Omega \to \R$ and $\Delta \Subset S$ that
\begin{align}\label{generalized-switching-lemma-identity}
        \int_{\Omega_\Delta}\int_{\Omega}c_{\Delta}(\omega, \xi_{\Delta})f(\omega)g(\xi_{\Delta}\omega_{\Delta^c})\mu(d\omega) \lambda_\Delta(d\xi_\Delta)
        = 
        \int_{\Omega_\Delta}\int_{\Omega}\hat{c}_{\Delta}(\omega, \xi_{\Delta})f(\xi_{\Delta}\omega_{\Delta^c})g(\omega)\mu(d\omega)\lambda_\Delta(d\xi_\Delta),  
    \end{align}
    where 
    \begin{align*}
       \hat{c}_\Delta(\eta, \xi_\Delta) = c_\Delta(\xi_\Delta \eta_{\Delta^c}, \eta_\Delta) \frac{\gamma_\Delta(\xi_\Delta  \lvert \eta_{\Delta^c})}{\gamma_\Delta(\eta_\Delta \lvert \eta_{\Delta^c})}. 
    \end{align*}
\end{lemma}
\noindent 
To keep the notation for conditional expectations in the upcoming proof simple, we will denote integration with respect to $\mu$ by $\bbE[\cdot]$.  
\begin{proof}
    As a first step, note that for fixed $\Delta \Subset S$ and $\xi_\Delta \in \Omega_\Delta$ the maps 
    \begin{align*}
        \Omega \ni \omega \mapsto g(\xi_{\Delta}\omega_{\Delta^c}) \in \R, \quad \Omega \ni \omega \mapsto f(\xi_{\Delta}\omega_{\Delta^c}) \in \R,
    \end{align*}
    are $\calF_{\Delta^c}$-measurable. 
    Therefore, we can use that $\gamma$ is the local conditional distribution of $\mu$ and the definition of the rates $\hat{c}$ to obtain the $\mu$-almost-sure identity 
    \begin{align*}
        \condexp{c_{\Delta}(\cdot, \xi_{\Delta})f(\cdot)g(\xi_{\Delta}\cdot_{\Delta^c})}{\calF_{\Delta^c}}(\omega) 
        &=
        g(\xi_{\Delta}\omega_{\Delta^c}) \condexp{c_{\Delta}(\cdot, \xi_{\Delta})f(\cdot)}{\calF_{\Delta^c}}(\omega)
        \nonumber \\\
        &=
        g(\xi_{\Delta}\omega_{\Delta^c}) 
        \int_{\Omega_\Delta}
            \gamma_{\Delta}(\zeta_{\Delta}|\omega_{\Delta^c})c_{\Delta}(\zeta_{\Delta}\omega_{\Delta^c}, \xi_{\Delta})f(\zeta_{\Delta}\omega_{\Delta^c})\lambda_\Delta(d\zeta_\Delta)
        \\\
        &= 
        g(\xi_{\Delta}\omega_{\Delta^c}) 
        \int_{\Omega_\Delta}
            \gamma_{\Delta}(\xi_{\Delta}|\omega_{\Delta^c})\hat{c}_{\Delta}(\xi_{\Delta}\omega_{\Delta^c}, \zeta_{\Delta})f(\zeta_{\Delta}\omega_{\Delta^c}) \lambda_\Delta(d\zeta_\Delta). 
        \nonumber 
    \end{align*}
    If we now integrate over $\xi_\Delta$, exchange the order of integration (via Fubini) and apply the same arguments in reverse – with $f$ taking the role of $g$ and vice versa – we get
    \begin{align*}
        \int_{\Omega_\Delta}
            \condexp{c_{\Delta}(\cdot, \xi_{\Delta})f(\cdot)g(\xi_{\Delta}\cdot_{\Delta^c})}{\calF_{\Delta^c}}(\eta) \lambda_\Delta(d\xi_\Delta)
        = 
        \int_{\Omega_\Delta}
            \condexp{\hat{c}_{\Delta}(\cdot, \zeta_{\Delta})f(\zeta_{\Delta}\cdot_{\Delta^c})g(\cdot)}{\calF_{\Delta^c}}(\eta)\lambda_\Delta(d\zeta_\Delta).
    \end{align*}
    By integrating both sides with respect to $\mu$, exchanging the order of integration, and applying the law of total expectation we obtain 
    \begin{align*}
       \int_{\Omega_\Delta}\int_{\Omega_\Delta}(\omega, \xi_{\Delta})f(\omega)g(\xi_{\Delta}\omega_{\Delta^c})\mu(d\omega) \lambda_\Delta(d\xi_\Delta)
        =
        \int_{\Omega_\Delta}\int_{\Omega}\hat{c}_{\Delta}(\omega, \zeta_{\Delta})f(\zeta_{\Delta}\omega_{\Delta^c})g(\omega)\mu(d\omega)\lambda_\Delta(d\zeta_\Delta),
    \end{align*}
    as desired.
\end{proof}

\subsection{Regularity of the time-reversal rates}
To make sure that $\hat{\mathscr{L}}$ is actually the generator of a well-defined interacting particle system we now show that the collection of transition measures $(\hat{c}_\Delta(\cdot, \cdot))_{\Delta \Subset S}$ satisfies the three conditions $\mathbf{(L1)}-\mathbf{(L3)}$. 
\begin{proposition}
    Let the rates of an interacting particle system with generator $\mathscr{L}$ satisfy $\mathbf{(R1)}-\mathbf{(R5)}$ and assume that $\mu$ is a time-stationary measure for $\mathscr{L}$ and such that $\mu$ is a Gibbs measure w.r.t.~a specification $\gamma$ that satisfies $\mathbf{(S1)}-\mathbf{(S4)}$.
    Then, the transition measures $(\hat{c}_\Delta(\cdot, d\xi_\Delta))_{\Delta \Subset \Z^d}$ with $\lambda_\Delta$-densities given by 
    \begin{align*}
        \hat{c}_\Delta(\eta, \xi_\Delta) = c_\Delta(\xi_\Delta \eta_{\Delta^c}, \eta_\Delta) \frac{\gamma_\Delta(\xi_\Delta  \lvert \eta_{\Delta^c})}{\gamma_\Delta(\eta_\Delta \lvert \eta_{\Delta^c})}
    \end{align*}
    satisfy the conditions $\mathbf{(L1)-(L3)}$. 
\end{proposition}
\begin{proof}
    \textit{Ad $\mathbf{(L1)}$:} This follows from the continuity assumptions $\mathbf{(R2)}$ and $\mathbf{(S2)}$, together with assumption $\mathbf{(S3)}$ and Lemma \ref{lemma:non-nullness-inequality}.
    \newline 
    \textit{Ad $\mathbf{(L2)}$:} 
    Note that for fixed $\Delta \Subset S$, $\xi_\Delta \in \Omega_\Delta$ and $\eta \in \Omega$ we have by Lemma \ref{lemma:non-nullness-inequality} and assumption $\mathbf{(R5)}$
    \begin{align*}
        \abs{c_\Delta(\eta, \xi_\Delta)} = \abs{c_\Delta(\xi_\Delta \eta_{\Delta^c},\eta_\Delta)\frac{\gamma_\Delta(\xi_\Delta \lvert \eta_{\Delta^c}}{\gamma_\Delta(\eta_\Delta \lvert \eta_{\Delta^c})}}
        \leq 
        \frac{1}{\delta}e^R c_\Delta(\xi_\Delta \eta_{\Delta^c}, \eta_\Delta). 
    \end{align*}
    So we get 
    \begin{align*}
        \sup_{\eta \in \Omega}\hat{c}_\Delta(\eta, \Omega_\Delta) 
        =
        \sup_{\eta \in \Omega}\int_{\Omega_\Delta}\hat{c}_\Delta(\eta, \xi_\Delta)\lambda_\Delta(d\xi_\Delta)
        &\leq \frac{1}{\delta}e^R \sup_{\eta \in \Omega}\int_{\Omega_\Delta}c_\Delta(\xi_\Delta \eta_{\Delta^c}, \eta_{\Delta})\lambda_\Delta(d\xi_\Delta)
        \\\
        &\leq 
       \frac{1}{\delta}e^R \sup_{\eta \in \Omega} \norm{c_\Delta(\eta, \cdot)}_\infty. 
    \end{align*}
    Therefore, assumption $(\mathbf{R3})$ implies that $\mathbf{(L2)}$ is also satisfied. 
    \newline 
    \textit{Ad $\mathbf{(L3)}$:} 
    Fix $\Delta \Subset S$, $y \in S$ and two configurations $\eta^1, \eta^2$ that only disagree at $y$. Then, for any $\xi_\Delta \in \Omega_\Delta$ we have 
    \begin{align*}
        \abs{\hat{c}_\Delta(\eta^1, \xi_\Delta) -\hat{c}_\Delta(\eta^2, \xi_\Delta)}
        &=
        \abs{
        c_\Delta(\xi_\Delta \eta^1_{\Delta^c}, \eta^1_\Delta)\frac{\gamma_\Delta(\xi_\Delta \lvert \eta^1_{\Delta^c})}{\gamma_\Delta(\eta^1_\Delta \lvert \eta^1_{\Delta^c})}
        -c_\Delta(\xi_\Delta\eta^2_{\Delta^c},\eta^2_\Delta)\frac{\gamma_\Delta(\xi_\Delta \lvert \eta^2_{\Delta^c})}{\gamma_\Delta(\eta^2_\Delta \lvert \eta^2_{\Delta^c})}
        }
        \\\
        &\leq 
        \abs{c_\Delta(\xi_\Delta \eta^1_{\Delta^c}, \eta^1_\Delta)}
        \abs{\frac{\gamma_\Delta(\xi_\Delta \lvert \eta^1_{\Delta^c})}{\gamma_\Delta(\eta^1_\Delta \lvert \eta^1_{\Delta^c})}
        -
        \frac{\gamma_\Delta(\xi_\Delta \lvert \eta^2_{\Delta^c})}{\gamma_\Delta(\eta^2_\Delta \lvert \eta^2_{\Delta^c})}
        }
        \\\
        &\qquad+
        \abs{\frac{\gamma_\Delta(\xi_\Delta \lvert \eta^2_{\Delta^c})}{\gamma_\Delta(\eta^2_\Delta \lvert \eta^2_{\Delta^c})}}
        \abs{c_\Delta(\xi_\Delta\eta^1_{\Delta^c}, \eta^1_{\Delta}) - c_\Delta(\xi_\Delta \eta^2_{\Delta^c}, \eta^2_{\Delta})}. 
    \end{align*}
    To estimate this further, we will have to make a case distinction over whether the site $y$ is contained in $\Delta$ or not. 
    If $y$ is contained in $\Delta$, then we can naively use Lemma \ref{lemma:non-nullness-inequality} and assumption $\mathbf{(R5)}$ to obtain the rough estimate 
    \begin{align*}
        \abs{\hat{c}_\Delta(\eta^1, \xi_\Delta) -\hat{c}_\Delta(\eta^2, \xi_\Delta)} 
        \leq 
        4 \frac{1}{\delta}e^R \sup_{\eta \in \Omega, \xi_\Delta \in \Omega_\Delta}\abs{c_\Delta(\eta, \xi_\Delta)} 
        \leq 
        \frac{4 e^R K(c)}{\delta}. 
    \end{align*}
    In the case where $y$ is not contained in $\Delta$, we can (and have to) be a bit more precise. Via the elementary algebraic rule 
    \begin{align*}
        ac - bd = \frac{1}{2}\left[(a-b)(c+d) + (a+b)(c-d)\right],
    \end{align*}
    and Lemma \ref{lemma:non-nullness-inequality} plus assumption $\mathbf{(R5)}$ one obtains 
    \begin{align*}
        &\abs{c_\Delta(\xi_\Delta \eta^1_{\Delta^c}, \eta^1_\Delta)}
        \abs{\frac{\gamma_\Delta(\xi_\Delta \lvert \eta^1_{\Delta^c})}{\gamma_\Delta(\eta^1_\Delta \lvert \eta^1_{\Delta^c})}
        -
        \frac{\gamma_\Delta(\xi_\Delta \lvert \eta^2_{\Delta^c})}{\gamma_\Delta(\eta^2_\Delta \lvert \eta^2_{\Delta^c})}
        }
        +
        \abs{\frac{\gamma_\Delta(\xi_\Delta \lvert \eta^2_{\Delta^c})}{\gamma_\Delta(\eta^2_\Delta \lvert \eta^2_{\Delta^c})}}
        \abs{c_\Delta(\xi_\Delta\eta^1_{\Delta^c}, \eta^1_{\Delta}) - c_\Delta(\xi_\Delta \eta^2_{\Delta^c}, \eta^2_{\Delta})}
        \\\
        &=
        \frac{1}{2}
    \abs{c_{\Delta}(\xi_{\Delta}\eta^1_{\Delta^c}, \eta^1_{\Delta})}
    \abs{\frac{1}{\gamma_{\Delta}(\eta^1_{\Delta}|\eta^1_{\Delta^c})\gamma_{\Delta}(\eta^2_{\Delta}|\eta^2_{\Delta^c})}}
    \abs{\gamma_{\Delta}(\xi_{\Delta}|\eta^1_{\Delta^c}) - \gamma_{\Delta}(\xi_{\Delta}|\eta^2_{\Delta^c})}
    \abs{\gamma_{\Delta}(\eta^1_{\Delta}|\eta^1_{\Delta^c})+\gamma_{\Delta}(\eta^2_{\Delta}|\eta^2_{\Delta^c})}
    \\\
    &\qquad+
    \abs{\frac{\gamma_{\Delta}(\xi_{\Delta}|\eta^2_{\Delta^c}) }{\gamma_{\Delta}(\eta^2_{\Delta}|\eta^2_{\Delta^c})}}
    \abs{c_{\Delta}(\xi_{\Delta}\eta^1_{\Delta^c}, \eta^1_{\Delta})-c_{\Delta}(\xi_{\Delta}\eta^2_{\Delta^c}, \eta^2_{\Delta})}
    \\\
    &\leq 
    \frac{1}{2\delta^2}e^{2R} K(c)K(\gamma)\abs{\gamma_{\Delta}(\xi_{\Delta}|\eta^1_{\Delta^c}) - \gamma_{\Delta}(\xi_{\Delta}|\eta^2_{\Delta^c})}
    +
    \frac{1}{\delta}e^R \abs{c_{\Delta}(\xi_{\Delta}\eta^1_{\Delta^c}, \eta^1_{\Delta})-c_{\Delta}(\xi_{\Delta}\eta^2_{\Delta^c}, \eta^2_{\Delta})}.
    \end{align*}
    Now, by integrating this pointwise difference of the densities over $\xi_\Delta$, we obtain via all of the other assumptions that
    \begin{align*}
        \sup_{x \in S}\sum_{\Delta \ni x}\sum_{y \neq x} \delta_y\hat{c}_\Delta < \infty.
    \end{align*}
    But this is precisely $\mathbf{(L3)}$ and the proof is finished. 
\end{proof}

With these two intermediate results at hand, we can now show that $\hat{\mathscr{L}}$ is indeed the generator of the time-reversal of $(\eta_t)_{t \geq 0}$ (w.r.t.~the time-stationary measure $\mu$).

\begin{proof}[Proof of Proposition \ref{proposition:time-reversal-ips}]
    It only remains to show that for all $f,g \in D(\Omega)$ we have 
    \begin{align*}
        \int_\Omega f(\omega) \mathscr{L}g(\omega) \mu(d\omega)
        =
        \int_\Omega \left(\hat{\mathscr{L}}f\right)(\omega)g(\omega)\mu(d\omega),
    \end{align*}
    since then the claimed time-reversal duality follows from Lemma \ref{lemma:time-reversal-generators}.

    For this, we first note that it suffices to show that the duality relation for the generators holds for all local functions $f,g: \Omega \to \R$. Indeed, if it holds for all pairs of local functions, we can then extend it to functions with bounded total oscillation by using the estimates from Lemma \ref{lemma:growth-estimate-triple-norm} and dominated convergence. Therefore, let $f,g$ be two local functions and let $\Lambda \Subset S$ be sufficiently large such that both $f$ and $g$ only depend on coordinates in $\Lambda$. 
    By first applying Lemma \ref{generalized-switching-lemma} and then using that $\mu$ is time-stationary with respect to the Markovian dynamics generated by $\mathscr{L}$, we see that 
    \begin{align*}
&\int_{\Omega}f(\omega)\mathscr{L}g(\omega)\mu(d\omega) 
        - 
        \int_{\Omega}\left(\hat{\mathscr{L}}f(\omega)\right)g(\omega)\mu(d\omega)
        \\\
        &=
        \sum_{\Delta \cap \Lambda \neq \emptyset}\int_{\Omega_\Delta}
            \int_{\Omega}c_{\Delta}(\omega, \xi_{\Delta})[f\cdot g(\xi_{\Delta}\omega_{\Delta^c})- f\cdot g(\omega)]\mu(d\omega)\lambda_\Delta(d\xi_\Delta)
        \nonumber = 
        \int_{\Omega}\mathscr{L} (f\cdot g) (\omega)\mu(d\omega) = 0,
    \end{align*}
    which finishes the proof.
\end{proof}

\section{Trajectorial decay of $\Phi$-entropies}\label{section:pathwise-phi-entropy-decay}
\noindent 
In this section we use the time-reversed process and a martingale argument to prove Theorem \ref{thm:trajectorial-decay-phi-entropy}. 
\subsection{The time-dependent martingale property}
The main technical tool will be the following lemma which can be seen as an extension of \cite[Lemma IV.20.12]{rogers_diffusions_2000} to our setting. 
\begin{lemma}\label{lemma:martingale-property}
Let $\mathscr{L}$ be the generator of an interacting particle system $(\eta(s))_{s \geq 0}$ such that its transition rates satisfy $(\mathbf{L1})-(\mathbf{L3})$ and let $\mu$ be a time-stationary measure w.r.t.~$\mathscr{L}$. Then, for all $f:[0,\infty) \times \Omega \to \R$ such that
\begin{enumerate}[i.]
    \item $f(\cdot, \eta) \in C^1([0,\infty))$ for all $\eta \in \Omega$ and
    \item for all $T > 0$ it holds that $$\sup_{0 \leq t \leq T}\vertiii{f(t,\cdot)} < \infty,$$ 
\end{enumerate} the process defined by 
\begin{align*}
    f(t, \eta(t)) - \int_0^t (\partial_s + \mathscr{L})f(s,\eta(s))ds 
\end{align*}
is a martingale w.r.t.~the filtration $\calG_t := \sigma(\eta(u) : 0 \leq u \leq t)$.
\end{lemma}
The proof of this lemma is not difficult but hard to find in the existing literature, we therefore give it in some detail. 
\begin{proof}
For functions $f$ as above, we define 
\begin{align*}
    M(s) := f(s,\eta(s)) - \int_0^s (\partial_u + \mathscr{L})f(u,\eta(u))du, \quad s \geq 0 
\end{align*}
and denote by $(P_t)_{t \geq 0}$ the Markov semigroup generated by $(\partial_s + \mathscr{L})$. Then, for $s<t$, the Markov property and the elementary identity 
\begin{align*}
    \frac{d}{dt}P_t = P_t(\partial_t + \mathscr{L}) = (\partial_t + \mathscr{L})P_t,
\end{align*}
give us 
\begin{align*}
    &\mathbb{E}\left[ 
        f(t, \eta(t)) - \int_0^t (\partial_u + \mathscr{L})f(u,\eta(u))du\Big\lvert \calG_s
    \right]
    \\\
    &=
    P_{t-s}f(s, \eta(s)) 
    -
    \int_0^s (\partial_u + \mathscr{L})f(u, \eta(u))du
    -
    \int_s^t P_{u-s}(\partial_u + \mathscr{L})f(s, \eta(s))du
    \\\
    &=
    P_{t-s}f(s, \eta(s)) 
    -
    \int_0^s (\partial_u + \mathscr{L})f(u, \eta(u))du
    -
    \int_s^t \frac{d}{du}P_{u-s}f(s, \eta(s))du
    \\\
    &= 
    f(s,\eta(s)) - \int_0^s (\partial_u + \mathscr{L})f(u,\eta(u))du.
\end{align*}
This shows that the process $(M(s))_{s \geq 0}$ is indeed a martingale. 
\end{proof}

This abstract tool now lets us establish the analogue of the first step in the case of a finite state space considered in Section~\ref{section:intro-trajectorial-idea}. 

\begin{proposition}\label{propositon:likelihood-martingale}Let $(\Omega, \calA, \mathbb{P})$ be a probability space on which the interacting particle system $(\eta(s))_{s \geq 0}$ is defined.
Denote the generator of $(\eta(s))_{s \geq 0}$ by $\mathscr{L}$ and assume that the rates satisfy $\mathbf{(R1)}-\mathbf{(R5)}$ and assume that under $\bbP$ we have $\eta(0) \sim\mu$, where $\mu$ is a time-stationary measure for the corresponding Markov semigroup $(P_t)_{t \geq 0}$ on $C(\Omega)$ that is generated by $\mathscr{L}$ and that $\mu$ is a Gibbs measure w.r.t.~a specification $\gamma$ that satisfies $\mathbf{(S1)}-\mathbf{(S4)}$.
Then, for all $f \in D(\Omega)$ and $T > 0$, the process defined by 
\begin{align*}
    P_{T-s}f(\eta(T-s)), \quad 0 \leq s \leq T,
\end{align*}
is a $((\hat{\calG}_t)_{0 \leq t \leq T}, \mathbb{P})$-martingale, where $\hat{\calG}_t = \sigma(\eta(T-s): \ 0 \leq s \leq t)$. 
\end{proposition}

\begin{proof}
    Note that by Lemma \ref{lemma:growth-estimate-triple-norm} we can apply Lemma \ref{lemma:martingale-property} to the function 
    \begin{align*}
        [0,T] \times \Omega \ni (s, \eta) \mapsto P_{T-s} f(\eta). 
    \end{align*}
    But since we have by the chain rule
    \begin{align*}
        \partial_s P_{T-s}f = -\hat{\mathscr{L}}P_{T-s}f,
    \end{align*}
    the correction term cancels out and we obtain the claimed martingale property. 
\end{proof}
\subsection{Trajectorial decay of $\Phi$-entropies}
With this preliminary result in place, we can now come to the proof of our main result. 
\begin{proof}[Proof of Theorem \ref{thm:trajectorial-decay-phi-entropy}]
\textit{Submartingale property:} By Jensen's inequality and Proposition \ref{propositon:likelihood-martingale} we immediately see that the process $(L^{\Phi,f}(t))_{t \geq 0}$, as defined in (\ref{thm:definition-phi-likelihood-process}), is a submartingale. 
\newline 
\textit{Doob--Meyer decomposition:} Here we want to apply Lemma \ref{lemma:martingale-property} to the function 
\begin{align*}
    g: [0,T] \times \Omega \ni (s, \eta) \mapsto \Phi(P_{T-s}f) \in \R. 
\end{align*}
Via the chain rule we see that 
\begin{align*}
    \partial_s g(s, \cdot) = \partial_s \Phi(P_{T-s}f(\cdot)) = -\Phi'(P_{T-s}f(\cdot))\hat{\mathscr{L}}P_{T-s}f(\cdot). 
\end{align*}
Applying the generator $\hat{\mathscr{L}}$ for fixed $s \in [0,T]$ yields 
\begin{align*}
    \hat{\mathscr{L}}g(s, \eta) = \sum_{\Delta \Subset S}\int_{\Omega_\Delta}\hat{c}_\Delta(\eta, \xi_\Delta)\left[\Phi(P_{T-s}f(\xi_\Delta \eta_{\Delta^c}))-\Phi(P_{T-s}f(\eta)) \right]\lambda_\Delta(d\xi_\Delta). 
\end{align*}
By putting these two ingredients together and using the previously introduced notation for the Bregman $\Phi$-divergence we obtain
\begin{align*}
    (\partial_s + \hat{\mathscr{L}})g(s,\eta) = \sum_{\Delta \Subset S}\int_{\Omega_\Delta}\hat{c}_\Delta(\eta, \xi_\Delta)\phidiv(P_{T-s}f(\xi_\Delta \eta_{\Delta^c}) \lvert P_{T-s}f(\eta))\lambda_\Delta(d\xi_\Delta).
\end{align*}
The claimed Doob--Meyer decomposition (\ref{thm:doob-meyer-decomposition}) of the submartingale $L^{\Phi,f}$ now follows from Lemma \ref{lemma:martingale-property}. 
\end{proof}

\section*{Acknowledgements}
The authors acknowledge the financial support of the Leibniz Association within the Leibniz Junior Research Group on \textit{Probabilistic Methods for Dynamic Communication Networks} as part of the Leibniz Competition.

\appendix 
\section{The time-reversal of Markov processes in equilibrium}
In this section, we briefly summarize some properties of the time-reversal of a Markov process w.r.t. a time-stationary measure. These results are classical but not particularly easy to find in the literature, at least in this formulation. 

We start by making precise what we mean by time-reversal of a stochastic process. Recall that any time-stationary stochastic process $(X_t))_{t \geq 0}$ can be extended (in law) to a process $(X_t)_{ -\infty < t < \infty}$ via Kolmogorov's extension theorem.
\begin{definition}
Let $(X_t)_{t \geq 0}$ be a time-stationary stochastic process. We call a process $(Y_t)_{t \geq 0}$ the \emph{time-reversal} of $X$ if 
\begin{align*}
    \text{\upshape Law}((X_{-t})_{t \geq 0}) = \text{\upshape Law}((Y_t)_{t \geq 0}). 
\end{align*}
\end{definition}
The intuition behind this definition is that forward in time the process $Y$ looks like the process $X$ run backwards. For Markov processes this notion can be characterized in terms of the semigroups and generators as follows. 
\begin{proposition}
Let $X= (X_t)_{t \geq 0}$ and $Y = (Y_t)_{t \geq 0}$ be Markov processes on a compact topological space $E$ with associated semigroups $(T_t)_{t \geq 0}$ and $(S_t)_{t \geq 0}$. Assume that $X$ has a time-stationary measure $\nu$ and for all $f,g \in C(E)$ we have 
\begin{align}\label{eqn:duality-semigroups}
    \int_E (T_t f)g d\nu = \int_E f (S_tg)d\nu.
\end{align}
Then $\nu$ is also time-stationary for $Y$ and $Y$ is the time-reversal of $X$ (w.r.t.~$\nu$). 
\end{proposition}

\begin{proof}
    By the duality relation (\ref{eqn:duality-semigroups}), $\nu$ is also time-stationary for $Y$ and we can extend it to a process $(Y_t)_{-\infty < t < \infty}$ with $Y_0 \sim \nu$. To show that $Y$ is the time-reversal of $X$ (w.r.t.~$\nu$) it suffices to show that for arbitrary $n \in \N$ and times $t_0 < t_2 <\cdots < t_n < t_{n+1}$ and functions $f_1,\dots, f_n \in C(E)$ it holds that 
    \begin{align*}
        \mathbb{E}_{Y_0 \sim \nu}\left[f_1(Y_{t_1}) \cdots  f_n(Y_{t_n})\right]
        =
        \mathbb{E}_{X_0 \sim \nu}\left[ f_1(X_{-t_1}) \cdots  f_n(X_{-t_n})\right]. 
    \end{align*}
    We will do this as follows. First, we introduce the notation $f_0 = f_{n+1} \equiv 1$ and define functions $g_0,\dots, g_{n+1}$ and  $h_0,\dots, h_{n+1}$ by $g_0 = h_{n+1} \equiv 1$,  
    \begin{align*}
        g_l(x) := \mathbb{E}_x\left[\prod_{k=0}^{l-1}f_k(X_{t_l - t_k})\right], \quad 1 \leq l \leq n+1
    \end{align*}
    and 
    \begin{align*}
        h_l(y) := \mathbb{E}_y\left[\prod_{k=l+1}^{n+1}f_k(Y_{t_k - t_l}) \right], \quad 0 \leq l \leq n.
    \end{align*}
    With this notation it suffices to show that the quantity 
    \begin{align*}
        \alpha_l := \int_E g_l(x)f_l(x)h_l(x)\nu(dx), \quad 0 \leq l \leq n+1,
    \end{align*}
    does not depend on $l$. Indeed, by stationarity of $X$ and $Y$, this will then imply that 
    \begin{align*}
        \mathbb{E}_{Y_0 \sim \nu}\left[f_1(Y_{t_1}) \cdots  f_n(Y_{t_n})\right]
        =
        \int_E h_0 d\nu 
        =
        \alpha_0
        =
        \alpha_{n+1}
        =
        \int_E g_{n+1}d\nu 
        = 
        \mathbb{E}_{X_0 \sim \nu}\left[ f_1(X_{-t_1}) \cdots  f_n(X_{-t_n})\right],
    \end{align*}
    exactly as we wanted. In order to show that $\alpha_l$ does not depend on $l$, we use the duality relation as follows. For $0 \leq l \leq n$ it holds that 
    \begin{align*}
        \alpha_l 
        &= \int_E g_l f_l h_l d\nu
        \\\
        &= \int_E g_l(x)f_l(x)T_{t_{l+1}-t_l}[f_{l+1}h_{l+1}](x)\nu(dx)
        \\\
        &= 
        \int_E S_{t_{l+1}-t_l}[g_l f_l](x) f_{l+1}(x)h_{l+1}(x)\nu(dx)
        \\\
        &=
        \int_E g_{l+1}f_{l+1}h_{l+1}d\nu 
        =
        \alpha_{l+1}.
    \end{align*}
    Therefore, $\alpha_l$ does not depend on $l$ and the claim follows.
\end{proof}

The duality relation (\ref{eqn:duality-semigroups}) for the semigroups can also be verified by checking a similiar property on the level of generators. 
The main technical tool for establishing this relation between the generator of a Markov process and its semigroup will unsurprisingly be the celebrated Hille--Yosida theorem which we recall here. 

\begin{theorem}[Hille--Yosida]\label{theorem:hille-yosida}
    There is a one-to-one correspondence between Markov generators on $C(E)$ and Markov semigroups on $C(E)$. This correspondence is explicitly given by 
    \begin{enumerate}
        \item[i.] \begin{align*}
            \domL = \set{f \in C(E): \ \lim_{t \downarrow 0} \frac{S_t f - f}{t} \ \text{exists}}, \quad \text{and} \\\
            \mathscr{L}f := \lim_{t \downarrow 0} \frac{S_t f-f}{t}, \quad f \in \domL. 
        \end{align*}
        \item[ii.]
        \begin{align*}
            S_t f = \lim_{n \to \infty}\left(\text{Id}-\frac{t}{n}\mathscr{L} \right)^{-n}f, \quad f \in C(E), t \geq 0. 
        \end{align*}
    \end{enumerate}
    Moreover, 
    \begin{enumerate}
        \item[iii.] if $f \in \domL$, then $S_t f \in \domL$ for all $t \geq 0$ and $\frac{d}{dt}S_t f = \mathscr{L}S_t f = S_t \mathscr{L}f$ and
        \item[iv.] for $g \in C(E)$ and $\lambda \geq 0$, the solution to the resolvent equation $f-\lambda \mathscr{L}f = g$ is given by 
        \begin{align*}
            f = \int_{0}^{\infty}e^{-t}S_{\lambda t}g dt.
        \end{align*}
    \end{enumerate}
\end{theorem}

\begin{proof}
    See Theorem 1.2.6 and Theorem 4.2.2 in \cite{ethier_markov_1986}.
\end{proof}

With this at hand, we can now formulate the time-reversal duality for generators. 
\begin{lemma}\label{lemma:time-reversal-generators}
Let $(T_t)_{t \geq 0}$ and $(S_t)_{t \geq 0}$ be two Markov semigroups on $C(E)$ where $E$ is a compact topological space with time-stationary measure $\nu$. Let $(\mathscr{A}, \dom{\mathscr{A}})$ and $(\mathscr{B}, \dom{\mathscr{B}})$ be their generators. If for all $f \in \dom{\mathscr{A}}$ and $g \in \dom{\mathscr{B}}$ it holds that 
\begin{align}\label{eqn:duality-generators}
    \int_E \left(\mathscr{A}f\right) g \ d\nu = \int_E f \left( \mathscr{B}g \right) \ d\nu,
\end{align}
then $(\ref{eqn:duality-semigroups})$ holds. It suffices if $(\ref{eqn:duality-generators})$ holds for $f,g$ in a core of the respective generators. 
\end{lemma}
\begin{proof}
    The duality relation implies that for all $f \in \dom{\mathscr{A}}$, $g \in \dom{\mathscr{B}}$ and $\lambda \geq 0$ it holds that
    \begin{align*}
        \int_E f(g - \lambda \mathscr{L}g)\ d\nu = \int_E g(f - \lambda \hat{\mathscr{L}}f)\ d\nu. 
    \end{align*}
    Now we can replace $f$ by $(\text{Id}-\lambda \mathscr{L})^{-1}f$ and $g$ by $(\text{Id}-\lambda \hat{\mathscr{L}})^{-1}g$ to see that for all $f,g \in C(\Omega)$
    \begin{align*}
        \int_E[(\text{Id}-\lambda \mathscr{A})^{-1}f]g \ d\nu 
        =
        \int_E[(\text{Id}-\lambda \mathscr{B})^{-1}g]f \ d\nu. 
    \end{align*}
    By iteratively applying this equality we obtain that for all $n \in \N$ and $\lambda \geq 0$ it holds that
    \begin{align*}
        \int_E[(\text{Id}-\lambda \mathscr{A})^{-n}f]g \ d\nu 
        =
        \int_E[(\text{Id}-\lambda \mathscr{B})^{-n}g]f \ d\nu. 
    \end{align*}
    For fixed $t \geq 0$ we can replace $\lambda$ by $t/n$ and use Part $ii.$ of Theorem~\ref{theorem:hille-yosida} to see that 
    \begin{align*}
        \lim_{n \to \infty}\int_E\left[(\text{Id}-\frac{t}{n} \mathscr{A})^{-n}f\right]g \ d\nu 
        =
        \int_E g T_tf \ d\nu
        \quad \text{ and }\quad
        \lim_{n \to \infty}\int_E\left[(\text{Id}-\frac{t}{n} \mathscr{B})^{-n}g\right]f \ d\nu
        =
        \int_E f S_tg \ d\nu,
    \end{align*}
    as desired.
\end{proof}

To sum this up, in order to show that a stationary Markov process $Y$ is the time-reversal of a stationary Markov process $X$, it suffices to check that their generators satisfy the duality relation (\ref{eqn:duality-generators}). 

\bibliography{references}

\begin{thebibliography}{BMDG05}

\bibitem[BMDG05]{banerjee_clustering_2005}
Arindam Banerjee, Srujana Merugu, Inderjit~S. Dhillon, and Joydeep Ghosh.
\newblock Clustering with {Bregman} {Divergences}.
\newblock {\em Journal of Machine Learning Research}, 6(58):1705--1749, 2005.

\bibitem[Cha04]{chafai_entropies_2004}
Djalil Chafaï.
\newblock Entropies, convexity, and functional inequalities, {On}
  {$\Phi$}-entropies and {$\Phi$}-{Sobolev} inequalities.
\newblock {\em Kyoto Journal of Mathematics}, 44(2), January 2004.

\bibitem[EK86]{ethier_markov_1986}
Stewart~N. Ethier and Thomas~G. Kurtz, editors.
\newblock {\em Markov {Processes}}.
\newblock Wiley {Series} in {Probability} and {Statistics}. John Wiley \& Sons,
  Inc., Hoboken, NJ, USA, March 1986.

\bibitem[FJ16]{fontbona_trajectorial_2016}
Joaquin Fontbona and Benjamin Jourdain.
\newblock A trajectorial interpretation of the dissipations of entropy and
  {Fisher} information for stochastic differential equations.
\newblock {\em The Annals of Probability}, 44(1), January 2016.

\bibitem[FV17]{friedli_statistical_2017}
Sacha Friedli and Yvan Velenik.
\newblock {\em Statistical {Mechanics} of {Lattice} {Systems}: {A} {Concrete}
  {Mathematical} {Introduction}}.
\newblock Cambridge University Press, 1 edition, November 2017.

\bibitem[Geo11]{georgii_gibbs_2011}
Hans-Otto Georgii.
\newblock {\em Gibbs measures and phase transitions}.
\newblock Number~9 in De {Gruyter} studies in mathematics. De Gruyter, Berlin ;
  New York, 2nd ed edition, 2011.
\newblock OCLC: ocn706965527.

\bibitem[HS76]{holley_l2_1976}
Richard~A. Holley and Daniel~W. Stroock.
\newblock {$L^2$} theory for the stochastic {Ising} model.
\newblock {\em Zeitschrift für Wahrscheinlichkeitstheorie und Verwandte
  Gebiete}, 35(2):87--101, 1976.

\bibitem[JK22]{jahnel_dynamical_2022}
Benedikt Jahnel and Jonas Köppl.
\newblock Dynamical {Gibbs} {Variational} {Principles} for {Irreversible}
  {Interacting} {Particle} {Systems} with {Applications} to {Attractor}
  {Properties}, May 2022.
\newblock arXiv:2205.02738 [math].

\bibitem[KMS21]{karatzas_trajectorial_2021}
Ioannis Karatzas, Jan Maas, and Walter Schachermayer.
\newblock Trajectorial dissipation and gradient flow for the relative entropy
  in {Markov} chains.
\newblock {\em Communications in Information and Systems}, 21(4):481--536,
  2021.

\bibitem[KST20]{karatzas_trajectorial_2020}
Ioannis Karatzas, Walter Schachermayer, and Bertram Tschiderer.
\newblock Trajectorial {Otto} calculus, March 2020.
\newblock arXiv:1811.08686 [math].

\bibitem[KY22]{kim_trajectorial_2022}
Donghan Kim and Lane~Chun Yeung.
\newblock A trajectorial approach to entropy dissipation for degenerate
  parabolic equations.
\newblock {\em arXiv:2210.16158}, October 2022.

\bibitem[Lig05]{liggett_interacting_2005}
Thomas~M. Liggett.
\newblock {\em Interacting {Particle} {Systems}}.
\newblock Classics in {Mathematics}. Springer Berlin Heidelberg, Berlin,
  Heidelberg, 2005.

\bibitem[RW00]{rogers_diffusions_2000}
L.~C.~G. Rogers and David Williams.
\newblock {\em Diffusions, {Markov} {Processes}, and {Martingales}}.
\newblock Cambridge University Press, 2 edition, April 2000.

\bibitem[TCY23]{tschiderer_trajectorial_2023}
Bertram Tschiderer and Lane Chun~Yeung.
\newblock A trajectorial approach to relative entropy dissipation of
  {McKean}–{Vlasov} diffusions: {Gradient} flows and {HWBI} inequalities.
\newblock {\em Bernoulli}, 29(1), February 2023.

\end{thebibliography}
\bibliographystyle{alpha}
\end{document}